\documentclass{conm-p-l}
\usepackage{amssymb, amsmath, amsfonts}
\usepackage[all]{xy}
\usepackage{bm}
\usepackage{hyperref}
\def\genfd{{\bm k}}
\def\id{{\rm id}}
\long\def\nodo#1{{}}
\def\gg{\mathfrak{g}}

\def\PPartial#1{\frac{\partial}{\partial(\partial^{#1})}}
\def\Hom{\operatorname{Hom}}

\def\Der{\operatorname{Der}}

\def\End{\operatorname{End}}
\def\hx{\hat{x}}

\def\nxpoint{\refstepcounter{subsection}\makepoint{\thesubsection}}
\def\nxsubpoint{\refstepcounter{subsubsection}%
  \makepoint{\thesubsubsection}}
\def\refpoint#1{{\rm\textbf{\ref{#1}}}}
\def\makepoint#1{\medbreak\noindent{\bf #1. }}

\theoremstyle{plain}

\theoremstyle{definition}

\def\MR#1{} 

\begin{document}

\title{Twisted exterior derivative for universal enveloping algebras $\mathbf{I}$}
\author{Zoran \v{S}koda}
\address{
Department of Teacher's Education, University of Zadar, Franje Tudjmana 24,
HR-23000 Zadar, Croatia
}
\email{zskoda@unizd.hr}
\address{
Theoretical Physics Division,
Institute Rudjer Bo\v{s}kovi\'{c}, Bijeni\v{c}ka cesta~54, P.O.Box 180,
HR-10002 Zagreb, Croatia
}


\begin{abstract}
  Consider any representation $\bm\phi$ of a finite-dimensional
Lie algebra $\mathfrak{g}$ by derivations of the completed symmetric algebra
$\hat{S}(\mathfrak{g}^*)$ of its dual.
Consider the tensor product of $\hat{S}(\mathfrak{g}^*)$
and the exterior algebra $\Lambda(\mathfrak{g})$. 
We show that the representation $\bm\phi$ extends canonically
to the representation $\tilde{\bm\phi}$ of that tensor product
algebra. We construct an exterior derivative on that algebra,
giving rise to a
twisted version of the exterior differential calculus with the
universal enveloping algebra in the role of the coordinate algebra.
In this twisted version, the commutators between the noncommutative
differentials and coordinates are formal power series in partial
derivatives. The square of the corresponding exterior derivative
is zero like in the classical case, but the graded Leibniz rule is
deformed.
\end{abstract}


\maketitle
\tableofcontents
\section{Preliminaries and basic notation}

\nxpoint Viewing the universal enveloping algebras
of Lie algebras as noncommutative deformations
of symmetric algebras in our earlier article
with {\sc S. Meljanac}~\cite{scopr}
we have also deformed the (completed) Weyl algebra of differential
operators, found deformed analogues of partial derivatives,
and studied the deformations of Leibniz rules, all parametrized
by certain datum which comes in many disguises as
``orderings'', ``representations by (co)derivations'',
``realizations by vector fields'' and ``coalgebra isomorphisms
between $S(\gg)$ and $U(\gg)$''(cf. also~\cite{AC2,Petracci} and
\cite{ldWeyl}, Chap. 10). Some noncommutative deformations of this
kind are interesting for physical applications, for example the
$\kappa$-deformed Minkowski space~\cite{AC2,MS}. Physical picture should
eventually include full-fledged differential geometry on such spaces,
and field theories on them; in particular the gauge theories
based on connections on noncommutative fiber bundles.
Our main motivation is to perform some early steps
on the mathematical side of this complex programme.

\nxpoint In the present article, I will extend this picture
to consistently include the exterior calculus
in a canonical way, given the datum mentioned above.
After the appearance of the first arXiv version of this article,
more general non-canonical approaches
to a noncommutative differential calculus
were found for some special Lie algebras in \cite{MSasaForms}.
We find it remarkable that our canonical extension exists, and that
it stems from the unique extension of the smash product structure
between the universal enveloping algebra and the space of
so-called deformed derivatives to a bigger
smash product algebra which includes
also the exterior algebra $\Lambda(\gg)$.
This new observation is the main reason for writing this paper,
with a hope that this exterior differential calculus
will enable a variant of differential geometry where the
base space is represented by the (spectrum of) the enveloping algebra $U(\gg)$.

\nxpoint ({\em Universal enveloping algebras as Hopf algebras.})
We work over a fixed field $\genfd$ of characteristic $0$. The unadorned tensor product $\otimes=\otimes_\genfd$ is over $\genfd$.
Taking a universal enveloping algebra is a functor from Lie to associative $\genfd$-algebras; hence for any Lie algebra $\mathfrak{g}$,
the Lie algebra map $\mathfrak{g}\to\mathfrak{g}\oplus\mathfrak{g}$,
given by $x\mapsto x\oplus x$, induces an algebra map
$\Delta:U(\mathfrak{g})\to U(\mathfrak{g}\oplus\mathfrak{g})\cong U(\mathfrak{g})\otimes U(\mathfrak{g})$, the terminal map $\mathfrak{g}\mapsto 0$ induces
the map $\epsilon:U(\genfd)\to U(0)\cong\genfd$, and the antihomomorphism
$\gg\to\gg$, $x\mapsto -x$
induces an homomorphism $\gamma:U(\gg)\to U(\gg)^{\mathrm{op}}$.
$U(\gg)$ is a Hopf algebra with comultiplication $\Delta$, counit $\epsilon$ 
and antipode $\gamma$ (\cite{BourbakiAlg,montg,ldWeyl}).
We assume that the reader is familiar with the Sweedler notation~\cite{montg}
$\Delta(h) = \sum_i h_{(1)i}\otimes h_{(2)i} = \sum h_{(1)}\otimes h_{(2)}$,
with or without an explicit summation index $i$. 

\nxpoint \label{pt:Weyl} ({\em Notation on generators, duals; Weyl algebras.})
From now on, we fix an $n$-dimensional Lie $\genfd$-algebra $\gg$
with $\genfd$-basis $\hx_1,\ldots,\hx_n$. The basis elements
are identified with the generators of the universal enveloping
algebra $U(\gg)$. We introduce the structure constants $C^s_{ij}$
defined via commutation relations $[\hx_i,\hx_j] = C^s_{ij}\hx_s$.

To distinguish the canonical copy of $\gg$
embedded in the symmetric algebra $S(\gg)$ from the copy in $U(\gg)$,
we denote the corresponding
basis $x_1,\ldots, x_n$; this emphasize that $x_i$-s
commute in $S(\gg)$, while tha generators with hat,
$\hx_i$-s in $U(\gg)$, do not. The dual basis of
$\gg^*$ is denoted by $\partial^1,\ldots,\partial^n$.

Due to an antiisomorphism between the geometric picture with vector fields around unit element of a Lie group and the algebraic picture promoted here, it is natural, following~\cite{halg}, to introduce Weyl algebra $A_{n,\genfd}$ in (unusual) contravariant notation as the free associative $\genfd$-algebra generated by symbols $x_i,\partial^j$ whre $i,j-1,\ldots,n$ modulo the ideal generated by elements $\partial^i\partial^j-\partial^j\partial^i$, $x_i x_j-x_j x_i$ and $\partial^ix_j-x_j\partial^i-\delta^i_j$ where $\delta$ is the Kronecker symbol. Given a multiindex $I = (i_1,\ldots,i_n)\in\mathbf{N}^n_0$, denote $x_I := x_1^{i_1}\cdots x_n^{i_k}$, $\partial^I := (\partial^1)^{i_1}\cdots(\partial^n)^{i_n}$ and $|I| := \sum_{k=1}^n i_k$.
The elements of the form $x_I\partial^J$ where $I,J$ run over all multiindices form a basis of $A_{n,\genfd}$ and there is an increasing filtration on $A_{n,\genfd}$~(\cite{coutinho}) by the degree of differential operator which is, for an element $D=\sum_{I,J} a_{I,J} x_I\partial^J\in A_{n,\genfd}$, the maximal $k\in\mathbf{N}_0$ such that there is $J$ with $|J|=k$ and $a_{I,J}\neq 0$. Completion with respect to this filtration is a topological $\genfd$-algebra denoted $\hat{A}_{n,\genfd}$. 

\nxpoint \label{pt:smash}
If $A$ is an associative $\genfd$-algebra and $H$ a Hopf algebra,
then an action $\triangleright : H\otimes A\to A$,
or the equivalent representation $\rho: H\to\End_\genfd(A)$,
is a left {\bf Hopf action} (synonym: $(A,\triangleright)$
is a left Hopf $H$-module algebra~\cite{montg})
if $h\triangleright(a\cdot b) = \sum (h_{(1)}\triangleright a)(h_{(2)}\triangleright b)$ (equivalently, $\rho(h)(a\cdot b) = \sum\rho(h_{(1)})(a)\cdot\rho(h_{(2)})(b)$) and $h\triangleright 1 = \epsilon(h) 1$ for all $h\in H$, $a,b\in A$. Similarly, one defines right Hopf actions.
Given a left Hopf action $\rho:H\to\End_\genfd(A)$,
the corresponding {\bf smash product algebra} $A\sharp_{\rho} H$
is an associative algebra with underlying vector space $A\otimes H$
and multiplication given by the unique linear extension of the formula
$$
(a \sharp g)(b\sharp h) = \sum a \rho(g_{(1)})(b) \sharp g_{(2)} h,
\,\,\,\,\,\,\,a,b\in A,\, g,h\in G,
$$
where $a\sharp h\in A\sharp_{\rho}H$ denotes $a\otimes h$
within the smash product. Similarly, for a right Hopf action $\sigma:H\to\End_\genfd^{\mathrm{op}}(A)$, there is a smash product algebra $H\sharp_\sigma A$ with
underlying space $H\otimes A$ and multiplication
$$
(g\sharp a)(h\sharp b) = \sum g h_{(1)} \sharp \sigma(h_{(2)})(a) b.
$$
Sometimes, it is useful to present $H\sharp_\sigma A$ by generators and relations. There are inclusions of algebras
$H\hookrightarrow H\sharp_\sigma A \hookleftarrow A$ given by
$H\ni h\mapsto h\sharp 1_A$ and $A\ni a\mapsto 1_H\sharp a$;
their images generate $H\sharp_\sigma A$.
This means that $H\sharp_\sigma A$ is a quotient of the free
product $H\ast A$ (the coproduct in the category of associative algebras). It is not difficult to see that the only additional relations are the
commutation relations between elements in $A$ and elements in $H$.
These may be read from the multiplication rule. Hence, 
$H\sharp_\sigma A$ is the quotient of the free product $H\ast A$ by the ideal generated by all elements of the form
$$
a h - \sum h_{(1)} \sigma(h_{(2)})(a),\,\,\,\,\,\,\,a\in A, h\in H. 
$$
In the case of a left Hopf action $\rho$, the relations are $h a - \sum \rho(h_{(1)})(a) h_{(2)}$.

\nxpoint A choice of basis $x_1,\ldots,x_n$
in an $n$-dimensional $\genfd$-vector space $V$
(e.g. $V=\mathfrak{g}$) induces a manifest isomorphism $S(V)\cong\genfd[x_1,\ldots,x_n]$ and of its dual $S(V)^* \cong \hat{S}(V^*)$ to the completed power ring $\genfd[[\partial^1,\ldots,\partial^n]]$.

The symmetric algebra on the dual $S(V^*)$ has a completion
(by the filtration by the degree of polynomial) $\hat{S}(V^*)$.
Algebra $\hat{S}(V^*)$ is isomorphic as a vector space
to the linear (algebraic) dual $(S(V))^*$: the nondegenerate pairing inducing this isomorphism is given by $\langle x_I,\partial^J \rangle = |I|! \delta^J_I$. In the interpretation of formal power series in
partial derivatives as representing
a formal differential operator of infinite order,
the duality pairing between $\hat{S}(V^*)$ and $S(V)$ is the evaluation of the differential operator at $0$. 

The symmetric algebra $S(V)$ is a Hopf algebra (even if $V$ is not finite dimensional) in a unique way such that
$\Delta(v) = 1\otimes v + v\otimes 1$ for all $v\in V$ (equivalently, $V$ may be given the structure of an Abelian Lie algebra $\mathfrak{a}$ and the canonical algebra isomorphism $S(V)\cong U(\mathfrak{a})$ transports the same coalgebra structure).

\nxpoint Let $A$ be an associative (unital) $\genfd$-algebra. We denote by $\Der(A)\subset\Hom_\genfd(A,A)$ the space of $\genfd$-linear derivations $A\to A$. It is a Lie $\genfd$-algebra with respect to the usual commutator.
Denote the algebra of $\genfd$-endomorphisms of a
$\genfd$-vector space $V$ by $\End_\genfd(V)$.
It is straightfoward to check that any homomorphism of Lie algebra
$\bm\phi:\gg\to\Der(A)$,
has a unique extension to a left Hopf action 
$U(\gg)\to\End(A)$. Similarly, any antihomomorphism
$\gg\to\Der(A)$ has a unique extension to a
right Hopf action $U(\gg)\to\End^{\mathrm{op}}(A)$. 

\nxpoint \label{pt:abelian}
If $\mathfrak{a}$ is finite dimensional Abelian Lie algebra with basis $x_1,\ldots,x_n$,
there is a Lie representation $\bm\delta:\mathfrak{a}\to\Der(S(\mathfrak{a}^*))$ (and the completed variant $\bm\delta:\mathfrak{a}\to\Der{\hat{S}(\mathfrak{a}^*)}$) given on generators $\partial^1,\ldots,\partial^n$
of $S(\mathfrak{a})$
by $\bm\delta(a)(\partial^i)= -\langle a,\partial^i\rangle$
and, in particular, $\bm\delta(x_i)(\partial^j) = -\delta_i^j$.
If $V$ is the underlying vector space of $\mathfrak{a}$, this induces a
right Hopf action
$\bm\delta\circ\gamma:S(V)\to\End^{\mathrm{op}}(S(V^*))$,
a completed variant $S(V)\to\End^{\mathrm{op}}(\hat{S}(V^*))$
and smash product algebras $S(V)\sharp_{\bm\delta\circ\gamma}S(V^*)\hookrightarrow S(V)\sharp_{\bm\delta\circ\gamma}\hat{S}(V^*)$.
The correspondence $x_i\mapsto x_i$, $\partial^j\mapsto\partial^j$
extends to an isomorphism between
the smash product algebra and the Weyl algebra $A_{n,\genfd}$
from~\refpoint{pt:Weyl} (the same holds for the completions). An alternative viewpoint are the symplectic Weyl algebras, where instead of the pairing, where the basic datum is a symplectic form on a module over a commutative ring ($V^*\oplus V$ is the module in our case), may be found in~\cite{ldWeyl}, Chap.\ 8.

\nxpoint \label{pt:derivpowerser} (\cite{BourbakiAlg})
All derivations of formal (commutative) power series rings,
and of $\hat{S}(\gg^*)$ in particular,
are continuous in formal topology.
Consider a $\genfd$-derivation $P$ and a formal power series
$D = D(\partial^1,\ldots,\partial^n)$.
Then the chain rule 
$P(D) = \sum_{r=1}^n\frac{\partial D}{\partial(\partial^r)} P(\partial^r)$
holds. The partial derivatives of $D$ are defined algebraically,
in the sense of formal power series.
In particular, the elements $P(\partial^r)$ for $r=1,\ldots,n$
determine the derivation $P$.

\nxpoint If $\bm\phi:\mathfrak{g}\to\Der(\hat{S}(\gg^*))$
is a Lie algebra homomorphism, then $\bm\phi(\hat{x}_i)\bm\phi(\hat{x}_j)(\partial^k) - \bm\phi(\hat{x}_j)\bm\phi(\hat{x}_i)(\partial^k) = \bm\phi([\hat{x}_i,\hat{x}_j])$ in particular, hence the matrix
$$\phi = (\phi^j_i) := (\bm\phi(-\hat{x}_i)(\partial^j))$$
with entries in $\hat{S}(\gg^*)$
satisfies the system of formal differential equations
\begin{equation}\label{eq:phiderphi}
\phi^l_j \PPartial{l}(\phi^k_i) - \phi^l_i \PPartial{l}(\phi^k_j)
= C^s_{ij}\phi^k_s,\,\,\,\,\,i,j,k\in\{1,\ldots,n\}
\end{equation}
Conversely, by linearity,
this system is sufficient for $\bm\phi$ to be a Lie homomorphism.
Therefore, homomorphisms $\bm\phi:\mathfrak{g}\to\Der(\hat{S}(\gg^*))$,
left Hopf actions $\bm\phi:U(\mathfrak{g})\to\Der(\hat{S}(\gg^*))$
and matrices $(\phi^i_j)$ satisfying~(\ref{eq:phiderphi}) are in 1-1
correspondence. Moreover, $\bm\phi\circ\gamma$ is a right Hopf action.
Recall that $\gamma(x) = -x$ for $x\in\mathfrak{g}$, hence
$(\bm\phi\circ\gamma)(x) = \bm\phi(-x)$ and $\phi^i_j = (\bm\phi\circ\gamma)(\hat{x}_i)(\partial^j)$.

\nxpoint \label{pt:assumphi}
From now on, we assume in addition that the matrix
$\phi$ is invertible. Consider 
the projection $\hat{S}(\gg^*)\to \hat{S}(\gg^*)/\cup_{n>0}
S^{n}(\gg^*)\cong\genfd$, sometimes interpreted as
the 'evaluation at $0$ map' (or the counit for the canonical structure of
Hopf algebra on the symmetric algebra).
If the image of $\phi^i_j$ under this
projection is the Kronecker $\delta^i_j$, we say that $\phi$ is
close to the unit matrix and symbolically write $\phi^j_i =
\delta^j_i + O(\partial)$. In the remainder of the paper, we assume
that $\phi$ is invertible and close to the unit matrix.
Given $\bm\phi$, these conditions on $\phi$ do not depend
on the choice of basis of $\mathfrak{g}$ used
in defining $\phi$, hence they are conditions on $\bm\phi$ only. 

\nxpoint  The left Hopf action $\bm\phi : U(\gg)\to\End(\hat{S}(\gg^*))$ and the right Hopf action $\bm\phi\circ\gamma$ induce
the isomorphic smash product algebras
$\hat{S}(\gg^*)\sharp_{\bm\phi} U(\gg)\cong
U(\gg)\sharp_{\bm\phi\circ\gamma}\hat{S}(\gg^*)$
which we may call the ``$\bm\phi$-twisted Weyl algebra''.
Indeed, by the free product description from~\refpoint{pt:smash},
both are quotients of the free product $\hat{S}(\gg^*)\ast U(\gg) =
U(\gg)\ast\hat{S}(\gg^*)$ by the commutation relations between
the generators $\hat{x}_i$ and arbitrary elements in $\hat{S}(\gg^*)$.
Using $\Delta(\hat{x}_i) = \hat{x}_i\otimes 1 + 1\otimes\hat{x}_i$, we obtain 
$D\hat{x}_i = \hat{x}_i \bm\phi(1)(D) + \bm\phi(\hat{x}_i)(D) = \hat{x}_i D
-\bm\phi(\hat{x}_i)(D)$ in one case and $\hat{x}_i D = \bm\phi(-\hat{x}_i)(D)+D\hat{x}_i$ in another smash product. 
By linearity of $\bm\phi$, these are clearly the same relations (the underlying general reason is that $U(\gg)$ is cocommutative, see~\cite{halg}, the discussion following Definition\ 1). For the generators, when $D = \partial^j$, these relations read $\partial^j\hat{x}_i - \hat{x}_i\partial^j = \phi^j_i$.

\nxpoint \label{pt:realization} ({\it Realizations in Weyl algebras.})
For a finite dimensional $\mathfrak{g}$ and general $\bm\phi$,
the correspondence $\hat{x}_i\mapsto \sum_{k=1}^n x_k\phi^k_i$,
$\partial^j\mapsto\partial^j$ does not depend on choice of basis,
and extends uniquely to a homomorphism
of algebras
$U(\gg)\sharp_{\bm\phi\circ\gamma}\hat{S}(\gg^*)\to\hat{A}_{n,\genfd}$.
If the matrix $\phi = (\phi^j_i)$ is invertible, this homomorphism
is invertible, with inverse given by $\hat{x}_i\mapsto \sum_{j=1}^n\hat{x}_k(\phi^{-1})^k_i, \partial^j\mapsto\partial^j$. 

\nxpoint It is shown in~\cite{scopr} that,
under the assumptions in~\refpoint{pt:assumphi},
the datum $\bm\phi$ is equivalent
to specifying a coalgebra isomorphism $\xi = \xi_{\bm\phi} :
S(\gg)\to U(\gg)$ (example: the symmetrization map
$\xi_{\mathrm{exp}}$~(\ref{eq:sym}) induced by $\phi^{\mathrm{exp}}$
from~(\ref{eq:univfla}))
which equals the identity
when restricted to $\genfd\oplus\gg\subset S(\gg)$.
However, there is no explicit formula relating $\bm\phi$ and $\xi$
in general. The isomorphism $\xi$ enables us to transport the
linear operators from $S(\gg)$ to $U(\gg)$. Partial derivatives $\partial^i$ transport to the deformed derivatives
$$\hat\partial^i = (\partial^i\!\blacktriangleright) :=
\xi\circ\partial^i\circ\xi^{-1}: U(\mathfrak{g})\to U(\mathfrak{g}),$$
satisfying the deformed Leibniz
rules studied in~\cite{scopr}. This action of generators $\partial^i$ on $U(\gg)$ together with action of elements in $U(\gg)$ on $U(\gg)$ by multiplication,
extend naturally to an action $\blacktriangleright$ of entire 
$U(\gg)\sharp_{\bm\phi\circ\gamma}\hat{S}(\gg^*)$ on $U(\gg)$.
In this paper, we try to avoid working 
with $\xi$ directly, but we do need $\hat{\partial}$. Fortunately,
there is an alternative description~\refpoint{pt:blackaction} of $\blacktriangleright$ in terms of the smash product, generators and matrix~$\phi$.

\nxpoint \label{pt:realizdet}
The $U(\gg)^*$ is a topological Hopf algebra by
duality with $U(\gg)$. The transpose
\begin{equation}\label{eq:transposexi}
  \xi^T:U(\gg)^*\to S(\gg)^*
\end{equation}
to any coalgebra isomorphism $\xi:S(\gg)\to U(\gg)$ is an algebra isomorphism.
Hence, one can transport the topological coalgebra structure
along this morphism
and obtain a nontrivial topological Hopf algebra structure on $\hat{S}(\gg^*)\cong S(\gg)^*$, studied in the disguise of deformed Leibniz rules
for $\hat{S}(\gg^*)$ in \cite{scopr}. Instead, we may transport the algebra structure from $U(\gg)$ to $S(\gg)$ obtaining the {\bf star product} $f\star g = \xi^{-1}(\xi(f)\cdot_{U(\gg)}\xi(g))$. The coproduct on $\hat{S}(\gg^*)$ is then the dual to the star product. An infinite dimensional version of the Heisenberg double construction applied to $U(\gg)$ results in an algebra, which is shown in~\cite{heisd} to be isomorphic to $U(\gg)\sharp_{\bm\phi\circ\gamma}\hat{S}(\gg^*)$ whenever the associated matrix $\phi$ satisfies the conditions from~\refpoint{pt:assumphi}. Algebra $U(\gg)$ is filtered by finite dimensional components and the dual is a cofiltered coalgebra. Thesis~\cite{stojicPhD} has exhibited that the Heisenberg double $U(\gg)^*\sharp U(\gg)$ is rigorously defined and has a structure of an internal Hopf algebroid in the symmetric monoidal category of filtered-cofiltered vector spaces. A variant of this Hopf algebroid with some extra completions (providing an {\em ad hoc} working setup, less satisfactory from the categorical point of view) has been exhibited a bit earlier in~\cite{halg}, in a form of $U(\gg)\sharp_{\bm\phi\circ\gamma}\hat{S}(\gg^*)$ and using explicit calculations with a concrete choice of $\bm\phi$. In that choice, the coalgebra isomorphism $\xi=\xi_{\mathrm{exp}}:S(\gg)\to U(\gg)$ is the symmetrization (or coexponential) map~\cite{BourbakiAlg,ldWeyl,halg,Petracci}, where for any $\hat{z}_1,\ldots,\hat{z}_k\in\genfd$ (not necessarily basis elements)
\begin{equation}\label{eq:sym}
\xi_{\mathrm{exp}}(z_1\ldots z_k) = \frac{1}{k!}\sum_{\sigma\in\Sigma(k)}\hat{z}_{\sigma(1)}\cdots\hat{z}_{\sigma(k)}
\end{equation}
or, alternatively $z^l\mapsto\hat{z}^l$ for any $z\in\mathfrak{g}\subset S(\gg)$ and $l\in\mathbf{N}$. The corresponding $\phi^i_j$ is
\begin{equation}\label{eq:univfla}
(\phi^{\mathrm{exp}})^i_j = \sum_{N=0}^\infty \frac{(-1)^N B_N}{N!}(\mathcal{C}^N)^i_j,
\end{equation}
where $B_N$ are the Bernoulli numbers and $\mathcal{C}$ is a matrix of elements in $\hat{S}(\gg^*)$ given by $\mathcal{C}^i_j = C^i_{jk}\partial^k$.
This formula has a long history (and direct relations to standard notions in Lie theory, e.g. the linear part of Hausdorff series) is derived over an arbitrary ring of characteristic $0$ in~\cite{ldWeyl}. In the case of $\genfd=\mathbf{R}$, the details of the relation to a geometric realization via formal differential operators around the unit element on a Lie group is exhibited in~\cite{halg},\ {\bf 1.2}.

\nxpoint\label{pt:blackaction}
Under the assumptions~\refpoint{pt:assumphi}, there is a deformed Fock action of $U(\gg)\sharp_{\bm\phi\circ\gamma}\hat{S}(\gg^*)$ on the universal enveloping algebra
$$
\blacktriangleright : (U(\gg)\sharp_{\bm\phi\circ\gamma}\hat{S}(\gg^*))\otimes U(\gg) \to U(\gg),
$$
obtained in four steps~(\cite{halg}). First we embed the second tensor factor $U(\gg)$ a $U(\gg)\otimes\genfd \hookrightarrow U(\gg)\sharp_{\bm\phi\circ\gamma}\hat{S}(\gg^*)$, then multiply, then use the isomorphism with $U(\gg)\sharp_{\bm\phi}\hat{S}(\gg^*)$ and then project the second factor to $\genfd$ via the counit of the symmetric algebra (the 'evaluation at $0$ map' from~\refpoint{pt:assumphi}); the result is in $U(\gg)\otimes\genfd\cong U(\gg)$. This is the standard physics procedure of pushing the partial derivatives to the right across coordinates (this time noncommutative coordinates, $\hat{x}_i$) using commutation relations and when all partial derivatives are on the right and coordinates on the left, retaining only the summand without partial derivatives. This way $U(\gg)$ becomes a left module of $U(\gg)\sharp_{\bm\phi\circ\gamma}\hat{S}(\gg^*)$, the $\phi$-deformed Fock module. Element $|0\rangle = 1_{U(\gg)}$ is the deformed Fock vacuum $|0\rangle_\gg=|0\rangle$ and we denote the result of the action on the vacuum by $x\mapsto x|0\rangle$. Map $\xi$ is simply the inclusion of algebras $S(\gg)\hookrightarrow U(\gg)\sharp_{\bm\phi\circ\gamma}\hat{S}(\gg^*)$ extending $x_i\mapsto \sum_{k=1}^n\hat{x}_k(\phi^{-1})^k_i$, followed by the action on the deformed vacuum.

If, in the first tensor factor, we restrict $\blacktriangleright$ on $\hat{S}(\gg^*)$ this is essentially the harpoon action of the topological dual algebra; this restriction is a topological left Hopf action. If, in the first tensor factor, we restrict $\blacktriangleright$ on $U(\gg)$, action $\blacktriangleright$ is simply the multiplication in $U(\gg)$.

Assuming the presentation of the Weyl algebra $A_{n,\genfd}$ as the smash product $S(V^*)\sharp_{\bm\delta\circ\gamma}S(V)$ for $V = \genfd[x_1,\ldots,x_n]$ (see~\refpoint{pt:abelian}), the standard Fock module of $A_{n,\genfd}$ is a special case of the corresponding deformed Fock module, namely $S(V)=U(\mathfrak{a})$ and $\mathfrak{a}$ is the Abelian Lie algebra with underlying space $V$, and we use the deformed Fock space construction for this case, without completions. The action, which we now denote by $\triangleright$, clearly agrees with the standard description as the action of a differential operator~(\cite{coutinho}) and the Fock action is $|0\rangle = 1_{S(V)}$.
For $\hat{A}_{n,\genfd}$ we use completions of course. The linear map
$\xi^{-1}:U(\gg)\to S(\gg)$ can be described as embedding $U(\gg)\to\hat{A}_{n,\genfd}$, $\hat{x}_i\mapsto\hat{x}_i^\phi := \hat{x}_k\phi^k_i$ followed by the action on the standard Fock vacuum $1_{S(\gg)}$. 
In particular, we obtain the description of the operator $\hat{\partial^i}= \xi\circ\partial\circ\xi^{-1} = \partial^i\blacktriangleright$ on $U(\gg)$ by iteratively pushing partial derivatives to the right using the commutation relation 
\begin{equation}\label{eq:blaction}\partial^i\blacktriangleright(\hat{x}_j\hat{u}) = \hat{x}_j(\partial^i\blacktriangleright\hat{u}) + \phi^i_j\blacktriangleright\hat{u},\,,\,\,\,\,\,\hat{u}\in U(\gg)
\end{equation}
and, at the end of the inductive procedure, retaining the summand without partial derivatives. Of course, $\phi^i_j$ is an infinite series, but commuting with $\hat{x}_k$ can drop the degree of each monomial only by $1$ (unless the commutator vanishes), hence we retain at each step only the summands in $\phi$ of the degree at most the degree of the noncommutative polynomial to the right of it. See an example~\refpoint{pt:mexample}.

\nxpoint \label{pt:deeper}
(This paragraph are remarks for deeper understanding and is not used in the rest of the article.) Under the assumptions~\refpoint{pt:assumphi}, the realization map from~\refpoint{pt:realization} supplies the isomorphism between the smash product $U(\gg)\sharp_{\bm\phi\circ\gamma}\hat{S}(\gg^*)$ and the usual completed Weyl algebra $\hat{A}_{n,\genfd}$, so one may ask why introducing this smash product at all. The action $\blacktriangleright$ and the corresponding deformed space make a difference, along with constructions derived from it. This action is Hopf with respect to the topological Hopf algebra structure on $\hat{S}(\gg^*)$ coming from the duality with $U(\gg)$. The map $\bm\phi$ provides a relation between $\blacktriangleright$ and the 'harpoon' action of the formal dual $U(\gg)^*$ on $U(\gg)$. The structure of the topological Hopf algebra on $\hat{S}(\gg^*)\cong U(\gg)^*$ together with $U(\gg)$ considered as a Hopf module over it (actually more than that, a braided commutative Yetter-Drinfeld module making sense in a monoidal category of filtered-cofiltered vector spaces) together induce both a smash product in the appropriate category and an additional structure of an internal Hopf $U(\gg)$-algebroid on that internal smash product. This is shown in~\cite{stojicPhD}. Earlier, it was sketched in~\cite{heisd} that this internal smash product is the same as an associative algebra with our smash product above (where the Hopf algebra is $U(\gg)$ and the Hopf module $U(\gg)^*$, rather than the other way around). A version of the Hopf algebroid (with some extra completions and somewhat {\em ad hoc} axioms on completions) is also derived in~\cite{halg}. If we replace $\bm\phi$ by another choice of homomorphism, say $\bm\psi$, we observe
($\xi^T$-s are from~(\ref{eq:transposexi})) isomorphisms of internal Hopf algebroids
\begin{equation}\label{eq:halgisom}\xymatrix{
    U(\gg)\sharp_{\bm\phi\circ\gamma}\hat{S}(\gg^*)\ar[rr]^{\id_{U(\gg)}\otimes(\xi_{\bm\phi}^T)^{-1}}&&
U(\gg)\sharp U(\gg)^* \ar[rr]^{\id_{U(\gg)}\otimes\xi^T_{\bm\psi}}&& U(\gg)\sharp_{\bm\psi\circ\gamma}\hat{S}(\gg^*),
}\end{equation}
where the Heisenberg double smash product Hopf algebroid in the middle is constructed using only the canonical topological Hopf pairing between Hopf algebra $U(\gg)$ and its canonical dual (toplogical) Hopf algebra $U(\gg)^*$ (and filtered cofiltered structures in place~\cite{stojicPhD}). 
Actions $\blacktriangleright_{\bm\phi}$ and $\blacktriangleright_{\bm\psi}$ get interchanged along the composition of isomorphisms~(\ref{eq:halgisom}) tensored by the identity on the additional $U(\gg)$ factor. However, the operator $\hat\partial^i$ depends on $\bm\phi$ as the image of $\partial^i$ under the composition isomorphism~(\ref{eq:halgisom}) is not $\partial^i$ in general. 

\section{The twisted algebra of differential forms}

\nxpoint \label{lem:homcheck}
    {\bf Lemma.} {\it Let $\mathfrak{g}$ be a Lie algebra, $A$ an
  associative algebra and $\rho:\mathfrak{g}\to\Der(A)$ a linear map.

  (i) If $g,h\in\mathfrak{g}$ and $\rho(g)\rho(h) - \rho(h)\rho(g)-\rho([g,h])$ vanishes when applied on each of two elements $a,b\in A$ then it vanishes on their product $ab$.

  (ii) Suppose that there is a family of algebra generators
  $\{a_\lambda\}_{\lambda\in\Lambda}$ of $A$ and a 
  subset $S\subset\mathfrak{g}$ which spans $\mathfrak{g}$,
  and such that for all $g,h\in S$,  
  $$\rho(g)\rho(h)(a_\lambda) - \rho(h)\rho(g)(a_\lambda) = \rho([g,h])(a_\lambda),$$
  Then $\rho$ is a representation of $\mathfrak{g}$ on $A$ by derivations.
  
  (iii) Suppose $\hat{A}$ is a topological associative $\genfd$-algebra containing $A$ as a dense subalgebra. If for each $h\in H$ derivation $\rho(h)$ extends (automatically uniquely) to a continuous $\genfd$-linear map $\rho'(h)$ of $\hat{A}$, then $\rho'$ is a representation of $\genfd$ by continuous $\genfd$-derivations iff $\rho$ is a representation. 
}

\nxpoint \label{pt:exterior} ({\it Notation on exterior algebras.})
In our constructions, it will be useful to distinguish notationally
the generators of two distinct copies of the classical exterior
algebra $\Lambda(\gg)$: in the first the generators
will be denoted by $d\hx_1,\ldots,d\hx_n$
and in the second by $dx_1,\ldots,dx_n$
(the latter copy first appears in \refpoint{s:correspondence}
and will be
denoted $\Lambda_{\mathrm{cl}}(\gg)$).
Both bases correspond to $\hx_1,\ldots,\hx_n$
under $\gg\hookrightarrow \Lambda(\gg)$.
Recall the convention from the introduction:
$\phi^j_i := \bm\phi(-\hat{x}_i)(\partial^j)$.

\nxpoint\label{pt:mainthm}
    {\bf Main theorem.} {\it
Any Lie homomorphism $\bm\phi :\gg\to\Der(\hat{S}(\gg^*))$
satisfying the assumptions from~\refpoint{pt:assumphi}
uniquely extends to a Lie homomorphism $\tilde{\bm\phi}:\gg\to
\Der(\Lambda(\gg)\otimes\hat{S}(\gg^*))$ satisfying
\begin{equation}\label{eq:mainObservation}
\tilde{\bm\phi}(\hx_i)(d\hx_l) = -\sum_{k,r,s=1}^n d\hx_k (\phi^{-1})^k_s
\left(\PPartial{r}\phi^s_l\right)\phi^r_i.
\end{equation}
This extension does not depend on the choice of basis of $\gg$.
}

{\it Proof.} By the Leibniz rule, any $\genfd$-derivation of an algebra $A$ is determined by its values on the generators of $A$, hence uniqueness follows even without the requirement that $\tilde{\bm\phi}$ be a Lie homomorphism. For the existence of $\tilde{\bm\phi}$ as a $\genfd$-linear map,
since the values on $\hat{S}(\gg^*)$ and on $\gg$ are predetermined, it remains to be shown that the extension by the Leibniz rule to the entire $\Lambda(\gg)\otimes\hat{S}(\gg)$ is well defined.
The only new nontrivial relation is antisymmetry
$d\hx_r\wedge d\hx_s = - d\hx_s \wedge d\hx_r$. The Leibniz rule gives
$$
\tilde{\bm\phi}(\hx_i)(d\hx_r\wedge d\hx_s) =
d\hx_k\wedge d\hx_s (\phi^{-1})^k_a
\left(\PPartial{b} \phi^a_r\right)\phi^b_i +
d\hx_r\wedge d\hx_k (\phi^{-1})^k_a
\left(\PPartial{b} \phi^a_s\right)\phi^b_i,
$$
where the right-hand side is evidently antisymmetric under the exchange
$(r\leftrightarrow s)$. 

It remains to show that $\tilde{\bm\phi}$ is automatically a representation (Lie homomorphism). 
Every derivation of $\hat{S}(\gg^*)$ is continuous
(\cite{BourbakiAlg}) and satisfies the chain rule,
which also makes problem of extension of derivations
from $S(\gg^*)$ to $\hat{S}(\gg^*)$ trivial.
Similar statements hold for $\Lambda(\gg)\otimes\hat{S}(\gg^*)$.
Thus we can apply Lemma~\refpoint{lem:homcheck}~(ii) to $A = \Lambda(\gg)\otimes S(\gg^*)$ and (iii) to $\hat{A} = \Lambda(\gg)\otimes\hat{S}(\gg^*)$ to assert that $\tilde{\bm\phi}$ is a Lie homomorphism iff for all $i,j,l\in\{1,\ldots,n\}$,  
\begin{equation}\label{eq:phitilde}
\tilde{\bm\phi}(\hx_i)\tilde{\bm\phi}(\hx_j)(d\hat{x}_l)
-\tilde{\bm\phi}(\hx_j) \tilde{\bm\phi}(\hx_i)(d\hat{x}_l)
- \tilde{\bm\phi}([\hx_i,\hx_j])(d\hat{x}_l) = 0.
\end{equation}
Using the Leibniz rule for $\tilde{\bm\phi}(\hx_i)$ we calculate,
ommiting the summation sign when summing over repeated indices, 
$$\begin{array}{lcl}
  \tilde{\bm\phi}(\hx_i)\tilde{\bm\phi}(\hx_j)(d\hx_l)& = &
  -\tilde{\bm\phi}(\hx_i)(d\hx_k)(\phi^{-1})^k_s\frac{\partial\phi^s_l}{\partial(\partial^r)}\phi^r_j - d\hx_k\bm\phi(\hx_i)((\phi^{-1})^k_s)\frac{\partial\phi^s_l}{\partial(\partial^r)}\phi^r_j
  \\ && \,\,\,\,-d\hx_k(\phi^{-1})^k_s\,\bm\phi(\hx_i)\left(\frac{\partial\phi^s_l}{\partial(\partial^r)}\right)\phi^r_j
  -d\hx_k(\phi^{-1})^k_s\frac{\partial\phi^s_l}{\partial(\partial^r)}\bm\phi(\hx_i)(\phi^r_j)
  \end{array}
$$
We first show that the first two summands mutually cancel.
By direct substitution of~(\ref{eq:phitilde}), the first summand becomes $+d\hx_{k'}(\phi^{-1})^{k'}_p\frac{\partial\phi^p_k}{\partial(\partial^{r'})}\phi^{r'}_i(\phi^{-1})^k_s\frac{\partial\phi^s_l}{\partial(\partial^r)}\phi^r_j$. By the chain rule (see~\refpoint{pt:derivpowerser}), $\bm\phi(\hx_j)(D) = -\sum_{t=1}^n\frac{\partial D}{\partial(\partial^p)}\phi^p_j$. This for $D=(\phi^{-1})^k_s$, together with the formula for the derivative of the inverse matrix we find that
$\bm\phi(\hx_j)((\phi^{-1})^k_s) = (\phi^{-1})^k_{s'}\frac{\partial\phi^{s'}_{r'}}{\partial(\partial^p)}(\phi^{-1})^{r'}_s\phi^p_j$,
hence the second summand above is
$$
-d\hx_k(\phi^{-1})^k_{s'}\frac{\partial\phi^{s'}_{r'}}{\partial(\partial^p)}\phi^p_j(\phi^{-1})^{r'}_s\frac{\partial\phi^s_l}{\partial(\partial^r)}\phi^r_j.
$$
The cancelation of the first two terms may be now observed after appropriately renaming dummy indices. By the chain rule, the third summand is
$$
+d\hx_k(\phi^{-1})^k_s\frac{\partial^2\phi^s_l}{\partial(\partial^p)\partial(\partial^r)}\phi^p_i\phi^r_j
$$
and is clearly antisymmetric under exchange $(i\leftrightarrow j)$.
Therefore, in~(\ref{eq:phitilde}), it cancels with the third summand in
$\tilde{\bm\phi}(\hx_j)\tilde{\bm\phi}(\hx_i)(d\hx_l)$. It remains to
consider the contributions from the 4th summand, giving
$$\begin{array}{l}
\tilde{\bm\phi}(\hx_i)\tilde{\bm\phi}(\hx_j)(d\hx_l)-
\tilde{\bm\phi}(\hx_j)\tilde{\bm\phi}(\hx_i)(d\hx_l) = \\
\,\,\,\,\,\,\,\,\,= d\hx_k(\phi^{-1})^k_s\frac{\partial\phi^s_l}{\partial(\partial^r)} \left(\frac{\partial\phi^r_j}{\partial(\partial^p)}\phi^p_i - 
\frac{\partial\phi^r_i}{\partial(\partial^p)}\phi^p_j\right)\\
\,\,\,\,\,\,\,\,\,\stackrel{(\ref{eq:phiderphi})}= -d\hx_k(\phi^{-1})^k_s \frac{\partial\phi^s_l}{\partial(\partial^r)} C_{ij}^{k'} \phi_{k'}^r\\
\,\,\,\,\,\,\,\,\,= \tilde{\bm\phi}([\hx_i,\hx_j])(d\hx_l).
\end{array}
$$
The tensorial notation is suggestive for the covariance properties with respect to the choice of basis. Let the prime indices denote a new basis, $B$ be a matrix of change of basis of $\gg$ in the sense that $\hx_i= \hx'_{k'} B^{k'}_i$ and the for the dual basis $\partial^i= (B^{-1})^i_{k'}\partial^{k'}$. By linearity of the map $\bm\phi$ and of 
$\genfd$-derivation $\bm\phi(-\hx_{i'})$ then
$$\phi^j_i = \bm\phi(-\hx_i)(\partial^j) =
(B^{-1})^j_{j'}B^{i'}_i\bm\phi(-\hx_{i'})(\partial^{j'})
= (B^{-1})^j_{j'}B^{i'}_i\phi^{j'}_{i'}.$$
This gives also
$x_i = x_{k'} B^{k'}_j(B^{-1})^j_{j'} B^{i'}_i\phi^{j'}_i\phi^{j'}_{i'}
= x_{k'} B^{k'}_i$. For the differentials, the embedding $\gg\hookrightarrow\Lambda(\gg)$ forces basis change $d\hx_i=d\hx'_{k'} B^{k'}_i$ extended by the usual tensoriality of higher exterior powers.
For the left hand side of~(\ref{eq:mainObservation}), by linearity, $\bm\phi(\hx_i)(d\hx_l) = B_i^{i'}B_l^{l'}\bm\phi(\hx_{i'})(d\hx_{l'})$.  Thus all ingredients in~(\ref{eq:mainObservation}) behave tensorially just as expected from the position of indices and the tensorial form is enough for the conclusion. Explicitly, we substitute the above component changes to the ingredients of the right hand side of~(\ref{eq:mainObservation}) to obtain
$$d\hx_{k'} B^{k'}_k(\phi^{-1})^{k''}_{s'} B^{s'}_{s}(B^{-1})^k_{k''}\left(\PPartial{r'}\phi^{s''}_{l'}\right)B^{l'}_l (B^{-1})^s_{s''}\phi^{r'}_{i'}B^{i'}_i (B^{-1})^r_{r'},$$
which gives the same change of basis coefficients as for the left hand side, after the contractions of numerical matrices $B$ and $B^{-1}$ above are accounted for. 

\nxpoint\label{pt:extendedalg}
    {\bf Corollary.} {\it Any representation $\bm\phi:\gg\to\Der(\hat{S}(\gg^*)$ of a finite dimensional Lie algebra $\gg$ by derivations on completed symmetric algebra of its dual has
a canonical  extension to a Hopf action of the form
$$\tilde{\bm\phi} :U(\gg)\to\End_\genfd(\Lambda(\gg)\otimes\hat{S}(\gg^*))$$
and satisfying~(\ref{eq:mainObservation}). In particular, the smash products
$$
U(\gg)\sharp_{\gamma\circ\tilde{\bm\phi}}(\Lambda(\gg)\otimes\hat{S}(\gg))
\cong
(\Lambda(\gg)\otimes\hat{S}(\gg))\sharp_{\tilde{\bm\phi}}U(\gg),
$$
are well-defined.
}

\nxsubpoint These two smash product algebras are canonically isomorphic by
the cocommutativity of $U(\gg)$. We call any of the two
the {\bf extended algebra of $\bm\phi$-twisted differential forms}.
'Extended' is for the additional presence
of partial derivatives in the algebra. 
A different, less canonical,
recipe is needed if we want that the subspace of
'forms' without partial derivatives be a subalgebra
(cf.~\cite{MSasaForms}). This is planned to be discussed in 
a sequel work (part II).
Notice that the usage of the antipode $\gamma$
in the choice of the smash product action
is in agreement with the minus sign in the formula
$\phi^j_i = \bm\phi(-\hat{x}_i)(\partial^j)$.

\nxpoint  We would like to describe the extended algebra of $\bm\phi$-differential forms by generators and relations. However, this algebra involves completions, and in particular formal power series in $\partial$-s and we need to consider it as either topological algebra or algebra with complete cofiltration (as a vector space) and multiplication distributing over formal sums in each argument, see~\cite{halg}, Appendix A.2 and~\refpoint{sp:cofilt}.

\nxsubpoint \label{sp:cofilt} A {\bf cofiltration} on a vector space $V$ is an inverse sequence of epimorphism of its quotients $\ldots \to V_i\to V_{i-1}\to \ldots \to V_0$. Denote the quotient maps $\pi_i:V\to V_i$ and bonding epimorphisms $\pi_{i,i+k}:V_{i+k}\to V_i$, the identities $\pi_{i} = \pi_{i,i+k}\circ \pi_{i+k}$ and $\pi_{i,i+k+l}=\pi_{i,i+k}\circ\pi_{i+k,i+k+l}$ hold. It is useful for us to consider cofiltrations as dual construction to increasing (ascending) filtrations like on $U(\gg)$ which involve inclusions of filtered components. However, many expositions treat cofiltrations in terms of 'descending filtrations' of kernels $\ker\pi_0\supset \ldots\supset \ker\pi_i\supset\ker\pi_{i+1}\supset\ldots$ with $V_i = V/\ker{\pi_i}$ and $\pi_{i.i+k}$ induced by inclusions. These kernels form a basis of neighborhoods of $0$ in a topological approach. 
 
A thread is a sequence $(v_r)_{r\in\mathbf{N}_0}\in \Pi_r V_r$ such that $v_r = \pi_{r,r+k}(v_{r+k})$. The set of threads is a vector space (with a canonical cofiltration) $\lim_i V_i = \hat{V}$ equipped with canonical {\it completion} map $V\to\hat{V}$, $v\mapsto (\pi_r(v))_r$. If this map is an isomorphism of vector spaces, we say that the cofiltration is {\bf complete} and the vector space complete cofiltered. A possibly infinite expression $\sum_{\lambda\in\Lambda} v_\lambda$ of elements in $V$ is a {\bf formal sum} if for any $i\in\mathbf{N}_0$ there are only finitely many $\lambda$ such that $\pi_i(v_\lambda)= 0$. A formal sum may be viewed as a useful representative of an element of $\hat{V}$~\cite{halg,stojicPhD}; the completion $\hat{V}$ may be viewed as a set of equivalence classes of formal sums. A linear map $f:V\to W$ of $\genfd$-vector spaces equipped with complete cofiltrations {\bf distributes over formal sums} if for every formal sum $\sum_\lambda v_\lambda$, the expression $\sum_\lambda f(v_\lambda)$ is also formal sum and $f(\sum_\lambda v_\lambda) = \sum_\lambda f(v_\lambda)$. Consider category $\mathrm{cfalg}_\otimes$ whose objects are associative $\genfd$-algebras $A$ equipped with complete cofiltration and such that the multiplication $\mu: A\otimes A\to A$ distributes over formal sums in each argument (i.e. $\mu(a,-)$ and $\mu(-,a)$ distribute over formal sums for each $a\in A$); the maps are $\genfd$-linear maps distributing over formal sums. This is a bit weaker structure than that of complete cofiltered algebras where the multiplication extends by definition to the completed tensor product, $A\hat\otimes A\to A$; the difference is just like between the continuity in each argument and joint continuity. It is reflected by the subscript $\otimes$ in $\mathrm{cfalg}_\otimes$.

\nxsubpoint Consider now the free associative $\genfd$-algebra 
$$F = \genfd\langle \hx_1,\ldots,\hx_n,d\hx_1,\ldots,d\hx_n,\partial^1,\ldots,\partial^n\rangle$$ 
on 3n symbols  $d\hx_i$, $\hx_i$ and $\partial^i$, for $i =1,\ldots, n$
and a subalgebra $F_{xd}$ generated by $\hx_i$-s and $d\hx_i$-s only. 
We also consider the free product of associative 
algebras $F_\star = F_{xd}\star\genfd[[\partial^1,\ldots,\partial^n]]$ where we allow for formal power series in $\partial$-s which mutually commute, but do not commute with other generators. There are obvious maps from $F$ and from $F_\star$ to $\bm\phi$-twisted differential forms; the map from $F$ is not surjective due completions, while the map from $F_\star$ is surjective, but this algebra is canonically completed in a way which is incompatible with quotienting to twisted differential forms (the kernel would not be closed). 

Let us define a category $\mathcal{K}$ as follows. The objects are maps of associative algebras $j : F\to A$ where $A$ is an object of $\mathrm{cfalg}_\otimes$ (\refpoint{sp:cofilt}) such that 
\begin{itemize}
\item $j(u\partial^{i_1}\cdots\partial^{i_{k+1}})\in\ker\pi_k$ for each $k\in\mathbf{N}_0$, $u\in F_{xd}$ and $i_1,\ldots,i_{k+1}\in\{1,\ldots,n\}$. This is a statement about a cofiltration on $A$. Notice that $u$ multiplies partial derivative symbols from the left. 
\item The following equations 
hold after applying $j$ (which we omit in writing)
\begin{equation}\label{eq:commall}\begin{array}{ccc}
[\hx_i,\hx_j] = C^k_{ij}\hx_k,&
[\partial^i,\partial^j] = 0,&
[\partial^j,\hx_i] = \phi^j_i,
\\ %
\left[\partial^j,d\hx_i\right] = 0,& [d\hx_j,\hx_i]=
d\hx_s(\phi^{-1})^s_r\left(\PPartial{l}\phi^r_j\right)\phi^l_i,
& \{d\hx_i,d\hx_j\} = 0
\end{array}\end{equation} where the commutator $[,]$
and the anticommutator $\{,\}$
are in the sense of associative algebras. Intuitively, these are required relations in a completion of $j(F)\subset A$.
\end{itemize}

Morphisms in $\mathcal{K}$ from $j:F\to A$ to $j':F\to A'$ are morphisms $\alpha:A\to A'$ in $\mathrm{cfalg}_\otimes$ such that $j'=\alpha\circ j$. 

\nxsubpoint \label{p:genrel} {\bf Theorem.} {\it (i) The extended
 algebra of $\bm\phi$-twisted differential forms is the (codomain of)
the universal initial object in category $\mathcal{K}$.

(ii) Suppose, $\omega\in\Lambda(\gg)\otimes\hat{S}(\gg^*)$, $k\in\mathbf{N}$ and $\hat{g}_{1},\ldots,\hat{g}_{k}\in\gg$. Then
\begin{equation}
  (\tilde{\bm\phi}\circ\gamma)(\hat{g}_{1}\cdots\hat{g}_k)(\omega)
  = [[\ldots[\omega,\hat{g}_1],\ldots],\hat{g}_k].
\end{equation}

(iii) Degree of a twisted differential form is a well defined
nonnegative integer.
}
\vskip .017in

{\it Proof (sketch).} We use the description~\refpoint{pt:smash} of the smash product as free product $U(\gg)\star(\Lambda(\gg)\otimes\hat{S}(\gg^*))$ modulo all commutation relations of the form
\begin{equation}\label{eq:smashproof}
  \bm\phi(u)(\omega) = \sum\hat{u}_{(1)}\omega\gamma(\hat{u}_{(2)}),
\end{equation}
where $\hat{u}\in U(\gg)$ and $\omega\in\Lambda(\gg)\otimes\hat{S}(\gg^*)$. It is however sufficient to include the relations where
$\hat{u},\omega$ are within chosen sets of generators. The kernel for the map
from $F_\star$ to the free product $U(\gg)\star(\Lambda(\gg)\otimes\hat{S}(\gg^*))$ is easy to observe. The relations (\ref{eq:commall}) enable to write every element in $j(F)$ as formal sum of the terms with partials on the right, hence we can see cofiltration on the subspace of $A$ consisting of all such elements (for every $j$) using $j(u\partial^{i_1}\cdots\partial^{i_{k+1}})\in\ker\pi_k$. We leave to the reader to show using these facts that the initial object $j:F\to A$  lifts to a surjection $F_\star\to A$ and, using this, to conclude the universal property for $F\to U(\gg)\sharp_{\gamma\circ\tilde\phi}(\Lambda(\gg)\otimes\hat{S}(\gg))$.

(i) Regarding that
$\Delta(\hx_i)=1\otimes\hx_i+\hx_i\otimes 1$, $S\hx_i = -\hx_i$,
the calculation is easy:
substitute $\hx_i$ for $\hat{u}$,
then $d\hx_j$ or $\partial^j$ for $\omega$,
and use identities $\tilde{\bm\phi}(\hx_i)(\partial^j) = -\phi^j_i$
and~(\ref{eq:mainObservation}).

For (ii), notice that $\bm\phi\circ\gamma$ is an antiisomorphism, $\bm\phi(\gamma(\hat{g}_i))(\nu) = \bm\phi(-\hat{g}_i)(\nu) =
[\nu,\hat{g}_i]$ for all $\nu\in\Lambda(\gg)\otimes\hat{S}(\gg^*)$, and proceed by induction. 

\nxpoint \label{s:correspondence} {\bf Theorem.} {\it If $\phi =
(\phi^i_j)$ is invertible and close to the identity (\refpoint{pt:assumphi}), the correspondence
\begin{equation}\label{eq:corrextA}
\hx_i\mapsto x_i^\phi := \sum_j x_j \phi^j_i, \,\,\,\,\,\,\,
\partial^i\mapsto \partial^i,\,\,\,\,\,\,\,
d\hx_i\mapsto d\hx_i^\phi := \sum_j dx_j \phi^j_i
\end{equation}
extends uniquely to an isomorphism
\begin{equation}
U(\gg)\sharp_{\tilde{\bm\phi}\circ\gamma}
(\Lambda(\gg)\otimes\hat{S}(\gg^*))\longrightarrow
\Lambda_{\mathrm{cl}}(\gg)\otimes \hat{A}_{n,\genfd}.
\end{equation}
}

\nxsubpoint \label{sp:intermediate}
We here used the notation conventions on exterior algebras~\refpoint{pt:exterior}. By abuse of notation, we may however by $dx_i$ denote also the preimage of $dx_i\in\Lambda_{\mathrm{cl}}(\gg)\otimes\hat{A}_{n,\genfd}$ under the isomorphism~(\ref{eq:corrextA}), namely $dx_i := d\hx_k(\phi^{-1})^k_i$. Notice that $[dx_i,\hat{x}_j] = 0$ and $[dx_i,\partial^j]=0$. It follows that we have a monomorphism $\Lambda_{\mathrm{cl}}(\gg)\otimes(U(\gg)\sharp_{\bm\phi\circ\gamma}\hat{S}(\gg^*))\to U(\gg)\sharp_{\tilde{\bm\phi}\circ\gamma}(\Lambda(\gg)\otimes\hat{S}(\gg^*))$ which is in fact an isomorphism because the generators of the form $dx_i,\partial^j,\hx_l$ generate $U(\gg)\sharp_{\tilde{\bm\phi}\circ\gamma}(\Lambda(\gg)\otimes\hat{S}(\gg^*))$. In this context, we call $\Lambda_{\mathrm{cl}}(\gg)\otimes(U(\gg)\sharp_{\bm\phi\circ\gamma}\hat{S}(\gg^*))$ the {\em intermediate algebra}.

\nxsubpoint 
Sometimes, it is good to consider realizations of 
$\Lambda_{\mathrm{cl}}(\gg)\otimes (U(\gg)\sharp_{\bm\phi\circ\gamma}\hat{S}(\gg^*))$
making it isomorphic to $\Lambda_{\mathrm{cl}}(\gg)\otimes\hat{A}_{n,\genfd}$:
here the generators of $\Lambda(\gg)$ are
also $dx_i$ and $U(\gg)\sharp_{\bm\phi\circ\gamma}\hat{S}(\gg^*)$ is
realized by $\hx_i^\phi = \sum_k x_k\phi^k_j$ and $\partial^i$ is just $\partial^i$ from $\hat{A}_{n,\genfd}$. Hence the intermediate algebra and the isomorphic extended algebra of $\bm\phi$-twisted forms do not depend on $\bm\phi$ as algebras (more is true, along the lines in~\refpoint{pt:deeper}).

\nxsubpoint Notation for generators $\hx_i^\phi = \sum_k
x_k \phi^k_i$, $d\hx_i^\phi=\sum_k dx_k
\phi^k_i$ extends to polynomials. The composition
$U(\gg)\hookrightarrow U(\gg)\sharp \hat{S}(\gg^*)\cong
\hat{A}_{n,\genfd}$ is an algebra monomorphism which agrees with
$()^\phi : u\mapsto u^\phi$. We do not use notation $()^\phi$ for
derivatives, because our isomorphism sends $\partial^i$ to
$\partial^i$ and this does not depend on $\bm\phi$. Thus
$d\hx_i^\phi = \sum_k dx_k \phi^k_i\in
\Lambda(\gg)\otimes \hat{A}_{n,\genfd}$. If we commute
elements within the image $\Lambda(\gg)\otimes
\hat{A}_{n,\genfd}$, we get the same commutators as
in~(\ref{eq:commall}), but in a realization, e.g. $[d\hx^\phi_i, \hx_j^\phi] =
d\hx_s^\phi(\phi^{-1})^s_r\left(\PPartial{l}\phi^r_i\right)\phi^l_j$.
Thus the theorem~\refpoint{s:correspondence} may be interpreted as
a realization of the extended algebra of
$\bm\phi$-twisted differential
forms in terms of ordinary differential forms and partial
derivatives (allowing infinite series).

\section{Exterior derivative}

\nxpoint \label{pt:dphi}
    {\it Definition.} $\bm\phi$-{\bf twisted exterior derivative}
is the $\genfd$-linear map given by
$$
\hat{d} = \hat{d}^{\bm\phi} := \sum_{k,j} d\hx_k (\phi^{-1})^k_j\hat\partial^j :
U(\gg)\sharp(\Lambda(\gg)\otimes \hat{S}(\gg^*))
\to U(\gg)\sharp(\Lambda(\gg)\otimes \hat{S}(\gg^*)),
$$
where $(\phi^{-1})^k_j$ acts by multiplication, while $\hat\partial^j=\partial^j\blacktriangleright$ acts on $U(\gg)$ tensor factor only, 
using $\bm\phi$-deformed Fock action $\blacktriangleright$ from~\refpoint{pt:blackaction}. We may write the exterior multiplication following $d\hx_k$ for emphasis; the multiplication is assumed.

\nxsubpoint Notice that $\partial^j$ and $\hat\partial^j$ do not commute, hence $(\phi^{-1})^k_j$ and $\hat\partial^j$ do not commute in general. By the discussion in~\refpoint{pt:deeper} not only the smash product, but the action $\blacktriangleright$ is also independent of $\bm\phi$ up to an isomorphism. However, the expression $\hat\partial^j$ depends on $\bm\phi$, see~\refpoint{pt:deeper}. Thus, $\hat{d}^{\bm\phi}$ may depend on $\bm\phi$ in general. 

\nxpoint\label{pt:mexample}
It is clear that this operator does not depend on the choice of basis.
If we realize $\hx_i$ as $\hx^\phi_i\in \hat{A}_{n,\genfd}$
then in this realization, the abstract $\hat{d}$ can also
be written as
$$
\hat{d} = \sum_{k,j} d\hx_k^\phi (\phi^{-1})^k_j \hat{\partial}^j
= \sum_j dx_j \hat{\partial}^j,
$$
once we interpret $\hat{\partial}^j$ properly in $\Lambda_{\mathrm{cl}}\otimes\hat{A}_{n,\genfd}$, which comes naturally if we work in terms of the intermediate algebra~\refpoint{sp:intermediate}, $\Lambda_{\mathrm{cl}}\otimes(U(\gg)\sharp_{\bm\phi}\hat{S}(\gg^*))$. This is different from the {\bf usual exterior derivative} $d =\sum_k dx_k \partial^k$ as $\partial^k$ acts on the factor $S(\gg)$ in  $\Lambda_{\mathrm{cl}}\otimes(S(\gg)\sharp_{\bm\delta\circ\gamma}\hat{S}(\gg^*))$, see~\refpoint{pt:Weyl}. Of course, via the inverse of the realization map, we could 
transport $d$ to our setup as well, but there is nothing essentially new.

\nxpoint {\bf Example.} Consider the series representation for $\phi^i_j$,
$$
\phi^i_j = \delta^i_j + A^i_{jr}\partial^{r} + \frac{1}{2}A^i_{jrs}\partial^r\partial^s + O(\partial^3).
$$
For the case $\phi=\phi^{\mathrm{exp}}$ (see~\refpoint{pt:realizdet}), formula~(\ref{eq:univfla}) gives $A^i_{jr} = \frac{1}{2}C^i_{jr}$ and $A^i_{jrs}=\frac{1}{6}C^i_{ps}C^p_{jr}$.

Using the procedure from the end of~\refpoint{pt:blackaction} and~(\ref{eq:blaction}),  we obtain
$\hat\partial^k(\hx_i)=\delta^k_i$ and
$$\begin{array}{lcl}
    \hat\partial^k(\hx_i\hx_j) &=& \hx_i\delta^k_j + \hx_j\delta^k_i + A^s_{ij}\\
    \hat\partial^k(\hx_i\hx_j\hx_l)&=& \hx_i\hx_j\delta^k_l+\hx_i\hx_l\delta^s_j + \hx_j\hx_j\delta^s_i+A^k_{jl}\hx_i+A^k_{ij}\hx_l+A^k_{il}\hx_j + A^k_{ir}A^r_{jl}+\frac{1}{2}A^k_{ijl}
\end{array}$$
Thus, $\hat{d}(\hx_i)= dx_i$, but
\begin{equation}\begin{array}{lcl}\label{eq:example}
    \hat{d}(\hx_i\hx_j) &=& (dx_j)\hx_i+(dx_i)\hx_j + A^k_{ij}dx_k\\
    \hat{d}(\hx_i\hx_j\hx_l)&=& (dx_l)\hx_i\hx_j+(dx_j)\hx_i\hx_l
    +(dx_i)\hx_j\hx_l
    \\ &&\,\,\,\,\,+dx_k(A^k_{jl}\hx_i+A^k_{ij}\hx_l+A^k_{il}\hx_j+A^k_{ir}A^r_{jl}+\frac{1}{2}A^k_{ijl})
\end{array}\end{equation}
Therefore, Leibniz rule is deformed for $\hat{d}$ in general.

In a $\bm\phi$-realization,~(\ref{eq:example}) differs from the application of the classical exterior derivative $d$ on the realization which simply reads
$$\begin{array}{l}
  d(\hx^\phi_i \hx^\phi_j) = (d\hx_i^\phi)\hx_j^\phi + \hx_i^\phi (d\hx_j^\phi)\\
  d(\hx^\phi_i \hx^\phi_j\hx^\phi_l) = (d\hx_i^\phi)\hx_j^\phi\hx_l^\phi + \hx_i^\phi (d\hx_j^\phi)\hx^\phi_l + \hx^\phi_i\hx^\phi_j d\hx^\phi_l.
\end{array}$$

\nxpoint \label{pt:d2} {\bf Theorem.} {\it
  (i) $\hat{d}^2 = 0$.

  (ii) $\hat{d}(d\hx_s\wedge\omega) = -d\hx_s\wedge\hat{d}\omega$ 
and $\hat{d}(d x_s\wedge\omega) = -d x_s\wedge\hat{d}\omega$;
}
\vskip .017in
  
{\it Proof.} (i) Regarding that $dx_m$ commute with elements in
$S(\gg^*)$,
\begin{equation}\label{eq:d2is0}
  \hat{d}^2 =  (dx_j\hat\partial^j)\wedge (dx_m\hat\partial^m)
  = dx_j\wedge dx_m \hat{\partial}^j \hat{\partial}^m
\end{equation}
is a contraction of the antisymmetric tensor $dx_j\wedge dx_m$ and the symmetric tensor $\hat{\partial}^j \hat{\partial}^m$, hence zero.
Note that $\hat\partial$-s do not commute with $\partial$-s nor $(\phi^{-1})^k_j$, and $d\hx_l$ does not commute with $\hat{\partial}^k$. Thus, $\hat{d}^2\neq d\hx_k\wedge d\hx_l (\phi^{-1})^k_j\hat{\partial}^j(\phi^{-1})^l_m \hat{\partial}^m$ and $ (\phi^{-1})^k_j\hat{\partial}^j(\phi^{-1})^l_m \hat{\partial}^m$ is not a symmetric tensor, unlike written in my earlier preprint.

(ii) The identities follow from the definition od $\hat{d}$, definition $dx_i = d\hx_k(\phi^{-1})^k_i$ and the antisymmetry of the product of generators in $\Lambda(\gg)$.

\nxpoint {\bf Proposition.}  {\it
  $\hat{d}$ preserves the subalgebra $\Lambda_U(\gg)$ generated by all $dx_i$ and $\hx_i$. There is an embedding
  $$\Lambda_{\mathrm{cl}}(\gg)\otimes U(\gg)\hookrightarrow\Lambda_{\mathrm{cl}}(\gg)\otimes (U(\gg)\sharp_{\bm\phi\circ\gamma}\hat{S}(\gg^*))\cong U(\gg)\sharp_{\tilde{\bm\phi}\circ\gamma}(\Lambda(\gg)\otimes\hat{S}(\gg^*))$$
whose image is $\Lambda_U(\gg)$.
}

\nxsubpoint We say that $\Lambda_U(\gg)$ is the algebra of {\bf $U$-twisted differential forms}. It is a subalgebra of the extended algebra of $\bm\phi$-twisted differential forms, or equivalenty of
the ``intermediate algebra''~\refpoint{sp:intermediate}.

 \nxpoint {\bf Proposition.} {\it The
usual Fock space action of $A_{n,\genfd}$ (and
$\hat{A}_{n,\genfd}$) on $S(\gg)$ extends to an action of
$\Lambda_{\mathrm{cl}}(\gg)\otimes(U(\gg)\sharp_{\bm\phi\circ\gamma}\hat{S}(\gg^*))\cong\Lambda_{\mathrm{cl}}(\gg)\otimes\hat{A}_{n,\genfd}$ on
$\Lambda_{\mathrm{cl}}(\gg)\otimes S(\gg)$ by multiplication in the first tensor
factor. The $\bm\phi$-deformed Fock space action $\blacktriangleright$ (see~\refpoint{pt:blackaction}) of $U(\gg)\sharp_{\bm\phi\circ\gamma}\hat{S}(\gg^*)$ on $U(\gg)$ extends to an action of the intermediate subalgebra
$\Lambda_{\mathrm{cl}}(\gg)\otimes
(U(\gg)\sharp_{\bm\phi\circ\gamma}\hat{S}(\gg^*))$ on its subalgebra
$\Lambda_{\mathrm{cl}}(\gg)\otimes U(\gg)$. This action restricts to the
multiplication on $\Lambda_{\mathrm{cl}}(\gg)\otimes U(\gg)$.
 }
 \vskip .02in
Keep in mind that elements from $\Lambda_{\mathrm{cl}}(\gg)$
commute with the elements in $\hat{S}(\gg^*)$ and $U(\gg)$,
while the elements in $U(\gg)$ do not commute with those in $\hat{S}(\gg^*)$.
The cyclic vectors for the two extended Fock spaces are still
$|0\rangle = 1_{S(\gg)}$ and $1_{U(\gg)}$.

\nxpoint \label{s:thvacuum} {\bf Theorem.} {\it
  Let $\omega\in \Lambda_U(\gg)$. Consider the action on usual vacuum
  $|0\rangle = 1_{\hat{S}(\gg)}$.

(i) Symbolically
$(\Lambda_U(\gg))^\phi|0\rangle = \Lambda_{\mathrm{cl}}(\gg)\otimes S(\gg)$.
More explicitly, the linear map sending $\omega\in \Lambda_U(\gg)$
to $\omega^\phi|0\rangle$ sends $\Lambda_U(\gg)$
into $\Lambda_{\mathrm{cl}}(\gg)\otimes S(\gg)$.
This action on vacuum in $\bm\phi$-realization
is an {\bf isomorphism of vector spaces}
(we asume $\bm\phi$ close to identity)
with inverse given by the $\bm\phi$-deformed Fock action on $|0\rangle_{U(\gg)}$, that is $\id_{\Lambda_{\mathrm{cl}}(\gg)}\otimes\xi_{\bm\phi} : \Lambda_{\mathrm{cl}}(\gg)\otimes S(\gg)\to \Lambda_{\mathrm{cl}}(\gg)\otimes U(\gg)\cong\Lambda_U(\gg)$.

(ii) $(\hat{d}\omega)|0\rangle = d(\omega|0\rangle)$,
where $d$ on the right hand side denotes the usual exterior derivative.

(iii) (Poincar\'e lemma) 
If  $\hat{d}(\omega) = 0$ then there is $\nu\in\Lambda_U(\gg)$ such that
$\hat{d}(\nu) =\omega$.
}
\vskip .017in

{\it Proof.} (i) and (ii) follow by direct check.
(iii) follows from the classical Poincar\'e lemma for $\Lambda_{\mathrm{cl}}(\gg)\otimes\hat{S}(\gg)$, the isomorphism in (i), property (ii) and Theorem~\refpoint{pt:d2},~(ii). 

\nxpoint \label{pt:tildestar} ({\it Star product on $\Lambda_U$.})
Using the isomorphism of vector spaces
$$\tilde\xi = \tilde\xi_{\bm\phi}=\id_{\Lambda(\gg)}\otimes\xi_{\bm\phi} = \Lambda_{\mathrm{cl}}(\gg)\otimes S(\gg)\to \Lambda_{\mathrm{cl}}(\gg)\otimes U(\gg)\cong\Lambda_U(\gg),$$
we can easily extend the star product $\star = \star_{\bm\phi}$ 
on $S(\gg)$ (see~\refpoint{pt:realizdet})
to an associative product $\wedge_\star=\wedge_{\star_{\bm\phi}}$ on
$\Lambda_{\mathrm{cl}}(\gg)\otimes S(\gg)$ given by
$\omega\wedge_\star\nu = \tilde\xi^{-1}(\tilde\xi(\omega)\wedge_{\Lambda_{\mathrm{cl}}(\gg)\otimes U(\gg)}\tilde\xi(\nu))$. 

\nxpoint \label{pt:tildestarprop} {\bf Proposition.} {\it
  If $\hat\omega,\hat\nu\in\Lambda_U(\gg)$ then
  $\hat{d}(\hat\omega\wedge\hat\nu) = \tilde\xi(d(\tilde\xi^{-1}(\hat\omega)\wedge_\star\tilde\xi^{-1}(\hat\nu))).$
}

{\it Proof.} For the deformed derivatives
$\hat\partial^i = (\partial^i\!\!\blacktriangleright)\, = \xi\circ\partial\circ\xi^{-1}$, hence it is immediate (\cite{scopr}) that $\hat\partial^i(\xi(f)\cdot_{U(\gg)}\xi(g)) = \xi(\partial^i(f\star g))$ for all $f,g\in S(\gg)$. By (bi)linearity it is sufficient to prove the statement when $\hat\omega = dx_{i_1}\wedge\cdots\wedge dx_{i_r}\xi(f)$ 
and $\hat\nu = dx_{j_1}\wedge\cdots\wedge dx_{j_s}\xi(g)$. Then by the definitions of $\hat{d}$ and of $\tilde\xi$ and $\tilde\star$, we obtain
$$\begin{array}{lcl}
  \hat{d}(\hat\omega\wedge\hat\nu) &=& dx_k \wedge dx_{i_1}\wedge\cdots\wedge dx_{i_r}\wedge dx_{j_1}\wedge\cdots\wedge dx_{j_s}(\xi\circ\partial^i\circ\xi^{-1})(\xi(f)\cdot_{U(\gg)}\xi(g))\\
  &=& \tilde\xi(dx_k \wedge dx_{i_1}\wedge\cdots\wedge dx_{i_r}\wedge dx_{j_1}\wedge\cdots\wedge dx_{j_s} \partial^k(f\star g))\\
  &=& \tilde\xi(dx_k\wedge \partial^k (\tilde\xi^{-1}(\tilde\omega)\wedge_\star\tilde\xi^{-1}(\tilde\nu)))\\
  &=& \tilde\xi(d(\tilde\xi^{-1}(\tilde\omega)\wedge_\star\tilde\xi^{-1}(\tilde\nu))).
\end{array}$$

\nxpoint ({\it Conclusion.}) Given a Lie algebra homomorphism~$\bm\phi:\mathfrak{g}\to\Der(\hat{S}(\gg^*))$ satisying assumptions~\refpoint{pt:assumphi} we exhibited in~\refpoint{pt:mainthm} its canonical extension to a Hopf action of $U(\gg)$ on $\Lambda(\gg)\otimes U(\gg)$ and a differential $\hat{d}_{\bm\phi}$ (defined in~\refpoint{pt:dphi}) on the induced smash product algebra~\refpoint{pt:extendedalg}, making it into a complex (by~\refpoint{pt:d2}, (i)), but not a differential graded algebra (the graded Leibniz rule is not satisfied). As an algebra, it does not depend on the choice of $\bm\phi$. This complex has a subcomplex $\Lambda_U(\gg)$ which, as an algebra, is isomorphic to $\Lambda_{\mathrm{cl}}(\gg)\otimes U(\gg)$ and thus deforms the differential graded algebra of polynomial differential forms $\Lambda(V)\otimes S(V)$ as an algebra and as a complex, but without Leibniz compatibility. The differential on $\Lambda_U(\gg)$ is related in~\refpoint{s:thvacuum} to the usual differential on polynomial differential forms via a $\bm\phi$-deformed Fock construction whose action $\blacktriangleright$ appears also in the study of Hopf algebroid structure of a Heisenberg double of $U(\gg)$. The star product on polynomials induced by the coalgebra isomorphism $\xi_{\bm\phi}:S(\gg)\to U(\gg)$ extends to the classical complex of polynomial differential forms: the product on $\Lambda_U(\gg)$ is transferred via Fock action (\refpoint{pt:tildestar}, \refpoint{pt:tildestarprop}). 

In a sequel paper, we address noncanonical modifications to our construction and their relation to alternative deformations, e.g.\ from~\cite{MSasaForms}, and to an approach using Drinfeld-Xu 2-cocycle twists of Heisenberg double Hopf algebroids.  

\nxpoint {\bf Acknowledgements.} The results have been written up as a preprint at IRB, Zagreb in Spring 2008. During a major revision in May 2020, I have been partly supported by the Croatian Science Foundation under the Project ``New Geometries for Gravity and Spacetime'' (IP-2018-01-7615).

\end{document}